**On the Octonion-like Associative Division Algebra**

Juhi Khalid, Martin Bouchard
School of Electrical Engineering and Computer Science, University of Ottawa,
800 King Edward, Ottawa K1N6N5 Canada

ABSTRACT

Using elementary linear algebra, this paper clarifies and proves some concepts about a recently introduced octonion-like associative division algebra over $\mathbb{R}$. This octonion-like algebra is actually the same as the split-biquaternion algebra, an even subalgebra of Clifford algebra $Cl_{4,0}(\mathbb{R})$. For two seminorms described in the paper, it is shown that the octonion-like algebra is a seminormed composition algebra over $\mathbb{R}$. Moreover, additional results related to the computation of inverse numbers in the octonion-like algebra are introduced, showing that the octonion-like algebra is a seminormed division algebra over $\mathbb{R}$, i.e., division by any number is possible as long as the two seminorms are non-zero. Additional results on normalization of octonion-like numbers and some involutions are also presented. The elementary linear algebra descriptions used in the paper also allow straightforward software implementations of the octonion-like algebra.

KEYWORDS: associative octonion-like algebra, octonion-like seminormed composition algebra, octonion-like seminormed division algebra, split-biquaternion algebra, hypercomplex algebra.

# I. INTRODUCTION

In [1], [2], an eight-dimensional octonion-like associative division algebra over $\mathbb{R}$ was recently introduced and presented as a normed algebra under geometric algebra criteria. There has been some controversy with this octonion-like algebra, as it was noted that a new normed division algebra over $\mathbb{R}$ could contradict Hurwitz's theorem [3]. The octonion-like algebra is actually the same as the split-biquaternion algebra, a sum of two quaternion algebras which is an even subalgebra of Clifford algebra $Cl_{4,0}(\mathbb{R})$ [4], [5]. The octonion-like algebra presented in [1], [2] is also found to be the same as the 1d-up approach to conformal geometric algebra presented in [6].

Using elementary linear algebra and two scalar seminorms, in this paper we explain that the associative octonion-like algebra over $\mathbb{R}$ is a seminormed composition algebra and it is also a seminormed division algebra over $\mathbb{R}$, i.e., division by any number is possible as long as the two seminorms are non-zero. Therefore, as an associative seminormed division algebra over $\mathbb{R}$, the octonion-like algebra does not contradict Hurwitz's theorem. It should be noted that the two scalar seminorms and the linear algebra used in this paper differ from the norm used in the original papers on the octonion-like algebra [1], [2], which were developed under geometric algebra criteria and where a multi-dimensional non-scalar norm was first defined and subsequently constrained to become a scalar.

The associativity property found in the octonion-like algebra (like all Clifford algebras but unlike octonions which are non-associative) can be an important criterion for the practical applicability of an algebra. A proof for the preservation of the two seminorms under multiplication of octonion-like numbers ($\mathbf{Z} = \mathbf{XY}$, $\|\mathbf{Z}\| = \|\mathbf{X}\|\|\mathbf{Y}\|$) is provided in the Appendix of this paper, using elementary linear algebra.

It will be shown in this paper that the computation of an inverse octonion-like number involves a matrix inversion, and that the case of a singular matrix only occurs when one of the two seminorms of an octonion-like number is zero. Additional results on normalization of octonion-like numbers and some involutions are also presented in the paper.

The explanation for the seminormed nature of the octonion-like algebra (and lack of contradiction with Hurwitz's theorem), the additional results related to the computation of inverse numbers in the octonion-like algebra (showing its seminormed division algebra over $\mathbb{R}$ status), the additional results on normalization of octonion-like numbers and some involutions, and the developments and proof using simple linear algebra are the contributions of this paper. The linear algebra descriptions also allow straightforward software implementations of the octonion-like algebra.

## II. OCTONION-LIKE ALGEBRA AND MULTIPLICATION RULES

Table I presents the multiplication rules for the octonion-like algebra. From the diagonal of the table, it can be seen that in an octonion-like number:

$$\mathbf{X} = X(0) + X(1)e_1 + X(2)e_2 + X(3)e_3 + X(4)e_4 + X(5)e_5 + X(6)e_6 + X(7)e_7 \qquad (1)$$

the imaginary dimensions correspond to coefficients $X(1)$ to $X(6)$, while coefficients $X(0)$ and $X(7)$ correspond to "scalar" or non-imaginary dimensions. We define $\mathbf{X}^*$ as the conjugate (or "reverse") of $\mathbf{X}$, changing the sign of the coefficients for the imaginary dimensions in $\mathbf{X}$, i.e., $X(1)$ to $X(6)$:

$$\mathbf{X}^* = X(0) - X(1)e_1 - X(2)e_2 - X(3)e_3 - X(4)e_4 - X(5)e_5 - X(6)e_6 + X(7)e_7 \qquad (2).$$

Table I Table of multiplication rules in the octonion-like algebra

| 1st factor (↓) and<br>2nd factor (→)<br>of octonion-like<br>number<br>multiplication | $e_0 = 1$ | $e_1$ | $e_2$ | $e_3$ | $e_4$ | $e_5$ | $e_6$ | $e_7$ |
|---|---|---|---|---|---|---|---|---|
| $e_0 = 1$ | 1 | $e_1$ | $e_2$ | $e_3$ | $e_4$ | $e_5$ | $e_6$ | $e_7$ |
| $e_1$ | $e_1$ | -1 | $e_3$ | $-e_2$ | $-e_5$ | $e_4$ | $e_7$ | $-e_6$ |
| $e_2$ | $e_2$ | $-e_3$ | -1 | $e_1$ | $e_6$ | $e_7$ | $-e_4$ | $-e_5$ |
| $e_3$ | $e_3$ | $e_2$ | $-e_1$ | -1 | $e_7$ | $-e_6$ | $e_5$ | $-e_4$ |
| $e_4$ | $e_4$ | $e_5$ | $-e_6$ | $e_7$ | -1 | $-e_1$ | $e_2$ | $-e_3$ |
| $e_5$ | $e_5$ | $-e_4$ | $e_7$ | $e_6$ | $e_1$ | -1 | $-e_3$ | $-e_2$ |
| $e_6$ | $e_6$ | $e_7$ | $e_4$ | $-e_5$ | $-e_2$ | $e_3$ | -1 | $-e_1$ |
| $e_7$ | $e_7$ | $-e_6$ | $-e_5$ | $-e_4$ | $-e_3$ | $-e_2$ | $-e_1$ | 1 |

Following these multiplication rules, the octonion-like multiplication can be detailed as:

$$\mathbf{Z} = \mathbf{XY} = \left( X(0) + X(1)e_1 + X(2)e_2 + X(3)e_3 + X(4)e_4 + X(5)e_5 + X(6)e_6 + X(7)e_7 \right)$$
$$\left( Y(0) + Y(1)e_1 + Y(2)e_2 + Y(3)e_3 + Y(4)e_4 + Y(5)e_5 + Y(6)e_6 + Y(7)e_7 \right)$$
$$\begin{aligned}
&= \left( X(0)Y(0) - X(1)Y(1) - X(2)Y(2) - X(3)Y(3) - X(4)Y(4) - X(5)Y(5) - X(6)Y(6) + X(7)Y(7) \right) \\
&+ \left( X(1)Y(0) + X(0)Y(1) - X(3)Y(2) + X(2)Y(3) + X(5)Y(4) - X(4)Y(5) - X(7)Y(6) - X(6)Y(7) \right) e_1 \\
&+ \left( X(2)Y(0) + X(3)Y(1) + X(0)Y(2) - X(1)Y(3) - X(6)Y(4) - X(7)Y(5) + X(4)Y(6) - X(5)Y(7) \right) e_2 \\
&+ \left( X(3)Y(0) - X(2)Y(1) + X(1)Y(2) + X(0)Y(3) - X(7)Y(4) + X(6)Y(5) - X(5)Y(6) - X(4)Y(7) \right) e_3 \\
&+ \left( X(4)Y(0) - X(5)Y(1) + X(6)Y(2) - X(7)Y(3) + X(0)Y(4) + X(1)Y(5) - X(2)Y(6) - X(3)Y(7) \right) e_4 \\
&+ \left( X(5)Y(0) + X(4)Y(1) - X(7)Y(2) - X(6)Y(3) - X(1)Y(4) + X(0)Y(5) + X(3)Y(6) - X(2)Y(7) \right) e_5 \\
&+ \left( X(6)Y(0) - X(7)Y(1) - X(4)Y(2) + X(5)Y(3) + X(2)Y(4) - X(3)Y(5) + X(0)Y(6) - X(1)Y(7) \right) e_6 \\
&+ \left( X(7)Y(0) + X(6)Y(1) + X(5)Y(2) + X(4)Y(3) + X(3)Y(4) + X(2)Y(5) + X(1)Y(6) + X(0)Y(7) \right) e_7 \\
&= Z(0) + Z(1)e_1 + Z(2)e_2 + Z(3)e_3 + Z(4)e_4 + Z(5)e_5 + Z(6)e_6 + Z(7)e_7
\end{aligned} \quad (3).$$

From (3), we see that for dimensions $e_1, e_2, e_3, e_4, e_5, e_6$ substituting $X(i)$ by $Y(i)$ and $Y(i)$ by $X(i)$ $0 \le i \le 7$ does not lead to the same value, therefore in general the octonion-like number product $\mathbf{Z} = \mathbf{XY} \ne \mathbf{YX}$ is non-commutative. Moreover, we can verify that:

$$(\mathbf{XY})^* = \mathbf{Y}^* \mathbf{X}^* \quad (4).$$

Also from (3), we can build a matrix product involving only real-valued numbers, where the matrix $\mathbf{M_X}$ includes coefficients from the left-side octonion-like number $\mathbf{X}$ in the product $\mathbf{Z} = \mathbf{XY}$:

$$\mathbf{Z}_r = \mathbf{M_X} \mathbf{Y}_r \quad (5),$$

with

$$\mathbf{Z}_r = \left[ Z(0)\ Z(1)\ Z(2)\ Z(3)\ Z(4)\ Z(5)\ Z(6)\ Z(7) \right]^T \quad (6)$$

$$\mathbf{M_X} = \begin{bmatrix}
+X(0) & -X(1) & -X(2) & -X(3) & -X(4) & -X(5) & -X(6) & +X(7) \\
+X(1) & +X(0) & -X(3) & +X(2) & +X(5) & -X(4) & -X(7) & -X(6) \\
+X(2) & +X(3) & +X(0) & -X(1) & -X(6) & -X(7) & +X(4) & -X(5) \\
+X(3) & -X(2) & +X(1) & +X(0) & -X(7) & +X(6) & -X(5) & -X(4) \\
+X(4) & -X(5) & +X(6) & -X(7) & +X(0) & +X(1) & -X(2) & -X(3) \\
+X(5) & +X(4) & -X(7) & -X(6) & -X(1) & +X(0) & +X(3) & -X(2) \\
+X(6) & -X(7) & -X(4) & +X(5) & +X(2) & -X(3) & +X(0) & -X(1) \\
+X(7) & +X(6) & +X(5) & +X(4) & +X(3) & +X(2) & +X(1) & +X(0)
\end{bmatrix} \quad (7)$$

$$\mathbf{Y}_r = \left[ Y(0)\ Y(1)\ Y(2)\ Y(3)\ Y(4)\ Y(5)\ Y(6)\ Y(7) \right]^T \quad (8).$$

The subscript $r$ in $\mathbf{Z}_r, \mathbf{Y}_r$ stands for "real-valued elements only". Alternatively, re-organizing the elements of the octonion-like product $\mathbf{Z} = \mathbf{XY}$ as:

$$\begin{aligned}
&\mathbf{Z} = \mathbf{XY} = \\
&\left( Y(0)X(0) - Y(1)X(1) - Y(2)X(2) - Y(3)X(3) - Y(4)X(4) - Y(5)X(5) - Y(6)X(6) + Y(7)X(7) \right) \\
&+\left( Y(1)X(0) + Y(0)X(1) + Y(3)X(2) - Y(2)X(3) - Y(5)X(4) + Y(4)X(5) - Y(7)X(6) - Y(6)X(7) \right) e_1 \\
&+\left( Y(2)X(0) - Y(3)X(1) + Y(0)X(2) + Y(1)X(3) + Y(6)X(4) - Y(7)X(5) - Y(4)X(6) - Y(5)X(7) \right) e_2 \\
&+\left( Y(3)X(0) + Y(2)X(1) - Y(1)X(2) + Y(0)X(3) - Y(7)X(4) - Y(6)X(5) + Y(5)X(6) - Y(4)X(7) \right) e_3 \quad (9) \\
&+\left( Y(4)X(0) + Y(5)X(1) - Y(6)X(2) - Y(7)X(3) + Y(0)X(4) - Y(1)X(5) + Y(2)X(6) - Y(3)X(7) \right) e_4 \\
&+\left( Y(5)X(0) - Y(4)X(1) - Y(7)X(2) + Y(6)X(3) + Y(1)X(4) + Y(0)X(5) - Y(3)X(6) - Y(2)X(7) \right) e_5 \\
&+\left( Y(6)X(0) - Y(7)X(1) + Y(4)X(2) - Y(5)X(3) - Y(2)X(4) + Y(3)X(5) + Y(0)X(6) - Y(1)X(7) \right) e_6 \\
&+\left( Y(7)X(0) + Y(6)X(1) + Y(5)X(2) + Y(4)X(3) + Y(3)X(4) + Y(2)X(5) + Y(1)X(6) + Y(0)X(7) \right) e_7
\end{aligned}$$

we can build a 2[nd] matrix form involving only real-valued numbers, where this time the matrix $\mathbf{P}_\mathbf{Y}$ includes coefficients from the right-side octonion-like number $\mathbf{Y}$ in the product $\mathbf{Z} = \mathbf{XY}$:

$$\mathbf{Z}_r = \mathbf{P}_\mathbf{Y} \mathbf{X}_r \quad (10),$$

with:

$$\mathbf{P}_\mathbf{Y} = \begin{bmatrix}
+Y(0) & -Y(1) & -Y(2) & -Y(3) & -Y(4) & -Y(5) & -Y(6) & +Y(7) \\
+Y(1) & +Y(0) & +Y(3) & -Y(2) & -Y(5) & +Y(4) & -Y(7) & -Y(6) \\
+Y(2) & -Y(3) & +Y(0) & +Y(1) & +Y(6) & -Y(7) & -Y(4) & -Y(5) \\
+Y(3) & +Y(2) & -Y(1) & +Y(0) & -Y(7) & -Y(6) & +Y(5) & -Y(4) \\
+Y(4) & +Y(5) & -Y(6) & -Y(7) & +Y(0) & -Y(1) & +Y(2) & -Y(3) \\
+Y(5) & -Y(4) & -Y(7) & +Y(6) & +Y(1) & +Y(0) & -Y(3) & -Y(2) \\
+Y(6) & -Y(7) & +Y(4) & -Y(5) & -Y(2) & +Y(3) & +Y(0) & -Y(1) \\
+Y(7) & +Y(6) & +Y(5) & +Y(4) & +Y(3) & +Y(2) & +Y(1) & +Y(0)
\end{bmatrix} \quad (11)$$

$$\mathbf{X}_r = \left[ X(0)\ X(1)\ X(2)\ X(3)\ X(4)\ X(5)\ X(6)\ X(7) \right]^T \quad (12).$$

As a special case of the previous detailed formulation for $\mathbf{Z} = \mathbf{XY}$, $\mathbf{Z} = \mathbf{XX}^*$ with the conjugate $\mathbf{X}^*$ as in (2) can be detailed as:

$$
\begin{aligned}
\mathbf{XX}^* &= \big(X(0) + X(1)e_1 + X(2)e_2 + X(3)e_3 + X(4)e_4 + X(5)e_5 + X(6)e_6 + X(7)e_7\big) \\
&\qquad \big(X(0) - X(1)e_1 - X(2)e_2 - X(3)e_3 - X(4)e_4 - X(5)e_5 - X(6)e_6 + X(7)e_7\big) \\
&= \big(X(0)X(0) + X(1)X(1) + X(2)X(2) + X(3)X(3) + X(4)X(4) + X(5)X(5) + X(6)X(6) + X(7)X(7)\big) \\
&+ \big(X(1)X(0) - X(0)X(1) + X(3)X(2) - X(2)X(3) - X(5)X(4) + X(4)X(5) + X(7)X(6) - X(6)X(7)\big)e_1 \\
&+ \big(X(2)X(0) - X(3)X(1) - X(0)X(2) + X(1)X(3) + X(6)X(4) + X(7)X(5) - X(4)X(6) - X(5)X(7)\big)e_2 \\
&+ \big(X(3)X(0) + X(2)X(1) - X(1)X(2) - X(0)X(3) + X(7)X(4) - X(6)X(5) + X(5)X(6) - X(4)X(7)\big)e_3 \\
&+ \big(X(4)X(0) + X(5)X(1) - X(6)X(2) + X(7)X(3) - X(0)X(4) - X(1)X(5) + X(2)X(6) - X(3)X(7)\big)e_4 \\
&+ \big(X(5)X(0) - X(4)X(1) + X(7)X(2) + X(6)X(3) + X(1)X(4) - X(0)X(5) - X(3)X(6) - X(2)X(7)\big)e_5 \\
&+ \big(X(6)X(0) + X(7)X(1) + X(4)X(2) - X(5)X(3) - X(2)X(4) + X(3)X(5) - X(0)X(6) - X(1)X(7)\big)e_6 \\
&+ \big(X(7)X(0) - X(6)X(1) - X(5)X(2) - X(4)X(3) - X(3)X(4) - X(2)X(5) - X(1)X(6) + X(0)X(7)\big)e_7
\end{aligned} \quad (13).
$$

It can be noted that all the imaginary dimensions in the octonion-like number $\mathbf{XX}^*$ have zero coefficients, so (13) reduces to:

$$
\begin{aligned}
\mathbf{XX}^* &= \big(X(0)X(0) + X(1)X(1) + X(2)X(2) + X(3)X(3) + X(4)X(4) + X(5)X(5) + X(6)X(6) + X(7)X(7)\big) \\
&+ \big(X(7)X(0) - X(6)X(1) - X(5)X(2) - X(4)X(3) - X(3)X(4) - X(2)X(5) - X(1)X(6) + X(0)X(7)\big)e_7 \\
&= \sum_{i=0}^{7} X^2(i) + \left(-2\sum_{i=1}^{3} X(i)X(7-i) + 2X(0)X(7)\right)e_7
\end{aligned} \quad (14).
$$

Unlike the general octonion-like product $\mathbf{XY}$, the product $\mathbf{XX}^*$ is commutative ( $\mathbf{XX}^* = \mathbf{X}^*\mathbf{X}$ and $(\mathbf{XX}^*)\mathbf{Y} = \mathbf{Y}(\mathbf{XX}^*)$), as can be observed from (9), (13) and (14).

The fact that the octonion-like product $\mathbf{XX}^*$ has non-zero coefficients only for dimensions $e_0 = 1$ and $e_7$ has a corresponding result for the matrix products $\mathbf{M}_\mathbf{X}^T\mathbf{M}_\mathbf{X}$ and $\mathbf{P}_\mathbf{X}^T\mathbf{P}_\mathbf{X}$: they are the sum of a diagonal matrix and an anti-diagonal matrix, whose non-zero elements are located at the positions of results with dimensions $e_0 = 1$ and $e_7$ in Table 1, and whose values are obtained from the $e_0 = 1$ and $e_7$ components in (14):

$$\mathbf{M}_{\mathbf{X}}^T\mathbf{M}_{\mathbf{X}} = \mathbf{P}_{\mathbf{X}}^T\mathbf{P}_{\mathbf{X}} = \begin{bmatrix} a & 0 & 0 & 0 & 0 & 0 & 0 & b \\ 0 & a & 0 & 0 & 0 & 0 & -b & 0 \\ 0 & 0 & a & 0 & 0 & -b & 0 & 0 \\ 0 & 0 & 0 & a & -b & 0 & 0 & 0 \\ 0 & 0 & 0 & -b & a & 0 & 0 & 0 \\ 0 & 0 & -b & 0 & 0 & a & 0 & 0 \\ 0 & -b & 0 & 0 & 0 & 0 & a & 0 \\ b & 0 & 0 & 0 & 0 & 0 & 0 & a \end{bmatrix} \qquad (15)$$

with

$$\begin{aligned} a &= \sum_{i=0}^{7} X^2(i) = \mathbf{X}_r^T\mathbf{X}_r \\ b &= -\sum_{i=1}^{6} X(i)X(7-i) + 2X(0)X(7) = -2\sum_{i=1}^{3} X(i)X(7-i) + 2X(0)X(7) = \mathbf{X}_r^T(\mathbf{J}\mathbf{X}_r) \end{aligned} \qquad (16)$$

and

$$\mathbf{J} = \begin{bmatrix} 0 & 0 & 0 & 0 & 0 & 0 & 0 & 1 \\ 0 & 0 & 0 & 0 & 0 & 0 & -1 & 0 \\ 0 & 0 & 0 & 0 & 0 & -1 & 0 & 0 \\ 0 & 0 & 0 & 0 & -1 & 0 & 0 & 0 \\ 0 & 0 & 0 & -1 & 0 & 0 & 0 & 0 \\ 0 & 0 & -1 & 0 & 0 & 0 & 0 & 0 \\ 0 & -1 & 0 & 0 & 0 & 0 & 0 & 0 \\ 1 & 0 & 0 & 0 & 0 & 0 & 0 & 0 \end{bmatrix} \qquad (17).$$

It can be noted that $\mathbf{J}\mathbf{X}_r$ is a flipped version of the conjugate octonion-like number $\mathbf{X}^*$ coefficients:

$$\mathbf{J}\mathbf{X}_r = \begin{bmatrix} X(7) & -X(6) & -X(5) & -X(4) & -X(3) & -X(2) & -X(1) & X(0) \end{bmatrix}^T \qquad (18).$$

We can then write $\mathbf{M}_{\mathbf{X}}^T\mathbf{M}_{\mathbf{X}}$ compactly:

$$\mathbf{M}_{\mathbf{X}}^T\mathbf{M}_{\mathbf{X}} = (\mathbf{X}_r^T\mathbf{X}_r)\mathbf{I} + (\mathbf{X}_r^T(\mathbf{J}\mathbf{X}_r))\mathbf{J} \qquad (19).$$

These results will be useful for the proof to be presented in the Appendix.

To complete this section, we also present a result using the matrix $\mathbf{M}_{\mathbf{X}}$ to compute the integer power of an octonion-like number:

$$\mathbf{X}^n = \underbrace{\mathbf{X}\cdots\mathbf{X}}_{n} \tag{20}$$

$$\mathbf{X}_r^n = \underbrace{\mathbf{M}_{\mathbf{X}}\cdots\mathbf{M}_{\mathbf{X}}}_{n-1}\mathbf{X}_r = \mathbf{M}_{\mathbf{X}}^{n-1}\mathbf{X}_r = \mathbf{S}\Lambda^{n-1}\mathbf{S}^{-1}\mathbf{X}_r \tag{21}$$

for $n \geq 1$ integer, and where $\mathbf{M}_{\mathbf{X}} = \mathbf{S}\Lambda\mathbf{S}^{-1}$ assumes that the diagonalization of the matrix $\mathbf{M}_{\mathbf{X}}$ is possible, with eigenvalues in the diagonal matrix $\Lambda$ and corresponding eigenvectors in the columns of $\mathbf{S}$.

## III. SEMINORMS FOR THE OCTONION-LIKE ALGEBRA

Consider the following non-negative scalar expression resembling the previous expression for the product $\mathbf{X}\mathbf{X}^*$ in (14). This expression is the square of two seminorms $\|\mathbf{X}\|$ :

$$\begin{aligned}\|\mathbf{X}\|^2 &= X(0)X(0) + X(1)X(1) + X(2)X(2) + X(3)X(3) + X(4)X(4) + X(5)X(5) + X(6)X(6) + X(7)X(7)\\ &\quad + \lambda X(7)X(0) - \lambda X(6)X(1) - \lambda X(5)X(2) - \lambda X(4)X(3) - \lambda X(3)X(4) - \lambda X(2)X(5) - \lambda X(1)X(6) + \lambda X(0)X(7)\\ &= \sum_{i=0}^{7} X^2(i) - 2\lambda\sum_{i=1}^{3} X(i)X(7-i) + 2\lambda X(0)X(7)\end{aligned} \tag{22}.$$

There is a seminorm with $\lambda = 1$ and a seminorm with $\lambda = -1$. Note that with this definition we have $\|\mathbf{X}^*\| = \|\mathbf{X}\|$ , with the conjugate $\mathbf{X}^*$ as in (2). Eq. (22) corresponds to seminorms and not norms because $\|\mathbf{X}\| = 0$ does not imply that $\mathbf{X} = 0$. These seminorms can also be expressed in compact vector form:

$$\|\mathbf{X}\|^2 = \mathbf{X}_r^T\mathbf{X}_r + \lambda\mathbf{X}_r^T(\mathbf{J}\mathbf{X}_r) \tag{23}.$$

It is proven in the Appendix that for these two seminorms the octonion-like algebra is a seminormed composition algebra, i.e., $\mathbf{Z} = \mathbf{X}\mathbf{Y}$ , $\|\mathbf{Z}\| = \|\mathbf{X}\|\|\mathbf{Y}\|$.

In the seminorms of (22) the term $\sum_{i=0}^{7} X^2(i)$ is always larger or equal to the amplitude of the other terms $\left|-2\sum_{i=1}^{3} X(i)X(7-i) + 2X(0)X(7)\right|$. If $X^2(i)$ $0 \le i \le 7$ either has some dominant values or is sparsely populated, then the term $\sum_{i=0}^{7} X^2(i)$ becomes dominant and the seminorms become similar to the Euclidian norm. The seminorms in (22) have a value of zero only when the following conditions are met:

- all the terms $|X(i)|$ which are non-zero have the same magnitude $|X(i)|$;
- for each non-zero term $X(i)$ there is a corresponding non-zero term $X(7-i)$ $0 \le i \le 7$
- all the non-zero products in $-2\sum_{i=1}^{3} X(i)X(7-i) + 2X(0)X(7)$ have the same polarity.

To illustrate this, consider the simple octonion-like product $(1-e_7)(1+e_7) = 1 - e_7e_7 = 1 - 1 = 0$. We must also have $\|(1-e_7)\|\|(1+e_7)\| = \|0\| = 0$, so either $\|(1-e_7)\|$ has a seminorm of zero (true if seminorm with $\lambda = 1$ is used) or $\|(1+e_7)\|$ has a seminorm of zero (true if seminorm with $\lambda = -1$ is used).

## IV. INVERSE NUMBERS IN THE OCTONION-LIKE ALGEBRA

The octonion-like algebra allows computation of inverses (and thus to perform divisions). However, computing the inverse of an octonion-like number (or performing divisions in general) requires a matrix inversion, as shown below. Consider an inverse $\mathbf{X}_i$ with the same form as the inverse used for complex number, quaternion or octonion algebras:

$$\mathbf{X}_i = \mathbf{X}^{-1} = \frac{\mathbf{X}^*}{\|\mathbf{X}\|^2} \quad (24).$$

Using the multiplication rules for octonion-like numbers and the previous result for $\mathbf{XX}^*$ in (14), for a right-side inverse we obtain:

$$\mathbf{X}\mathbf{X}_i = \frac{\mathbf{X}\mathbf{X}^*}{\left\|\mathbf{X}\right\|^2} = \frac{\sum_{i=0}^{7} X^2(i) + \left(-2\sum_{i=1}^{3} X(i)X(7-i) + 2X(0)X(7)\right) e_7}{\left\|\mathbf{X}\right\|^2} \qquad (25),$$

and with the definition of the seminorms for octonion-like numbers in (22), this becomes:

$$\mathbf{X}\mathbf{X}_i = \frac{\mathbf{X}\mathbf{X}^*}{\left\|\mathbf{X}\right\|^2} = \frac{\sum_{i=0}^{7} X^2(i) + \left(-2\sum_{i=1}^{3} X(i)X(7-i) + 2X(0)X(7)\right) e_7}{\sum_{i=0}^{7} X^2(i) - 2\lambda\sum_{i=1}^{3} X(i)X(7-i) + 2\lambda X(0)X(7)} \neq 1 \qquad (26).$$

We see that this is not the right-side inverse of an octonion-like number such that $\mathbf{X}\mathbf{X}_i = 1$, although $\left\|\mathbf{X}\right\|\left\|\mathbf{X}_i\right\| = \left\|\mathbf{X}\right\|\frac{\left\|\mathbf{X}^*\right\|}{\left\|\mathbf{X}\right\|^2} = \left\|\mathbf{X}\right\|\frac{\left\|\mathbf{X}\right\|}{\left\|\mathbf{X}\right\|^2} = 1$, i.e., it is a correct right-side inverse if only the seminorm is considered.

A right-side inverse can be obtained by writing the octonion-like product $\mathbf{X}\mathbf{X}_i = 1$ in the matrix form introduced earlier in (5):

$$\mathbf{M}_{\mathbf{X}}\mathbf{X}_{i,r} = \left[1\ 0\ 0\ 0\ 0\ 0\ 0\ 0\right]^T \qquad (27)$$

with

$$\mathbf{X}_{i,r} = \left[X_i(0)\ X_i(1)\ X_i(2)\ X_i(3)\ X_i(4)\ X_i(5)\ X_i(6)\ X_i(7)\right]^T \qquad (28).$$

Therefore, we have an 8 by 8 linear set of equations to solve (which can possibly be implemented more efficiently making use of the structure in $\mathbf{M}_{\mathbf{X}}$, but this is not investigated here). Or, if we explicitly use the matrix inverse $\mathbf{M}_{\mathbf{X}}^{-1}$, we have:

$$\mathbf{X}_{i,r} = \mathbf{M}_{\mathbf{X}}^{-1}\left[1\ 0\ 0\ 0\ 0\ 0\ 0\ 0\right]^T \qquad (29).$$

Similarly, the solution for a left-side inverse $\mathbf{X}_i\mathbf{X} = 1$ can be directly obtained by writing the octonion-like product $\mathbf{X}_i\mathbf{X} = 1$ with the matrix form introduced earlier in (10):

$$\mathbf{P}_{\mathbf{X}}\mathbf{X}_{i,r} = \left[1\ 0\ 0\ 0\ 0\ 0\ 0\ 0\right]^T \qquad (30)$$

$$\mathbf{X}_{i,r} = \mathbf{P}_{\mathbf{X}}^{-1}\left[1\ 0\ 0\ 0\ 0\ 0\ 0\ 0\right]^T \qquad (31).$$

However, the solution $\mathbf{X}_i$ for the right-side inverse $\mathbf{X}\mathbf{X}_i = 1$ is the same as for the left-side inverse $\mathbf{X}_i\mathbf{X} = 1$, as any product $\mathbf{XY} = \mathbf{YX}$ is commutative when the result only has $e_0$ and $e_7$ components in (9).

To evaluate when $\mathbf{M}_\mathbf{X}$ or $\mathbf{P}_\mathbf{X}$ can be singular, the eigenvalues of $\mathbf{M}_\mathbf{X}$ or $\mathbf{P}_\mathbf{X}$ are used:

$$
\begin{aligned}
\lambda_0,\lambda_1 &= X(0)+X(7) \\
&\quad + i\sqrt{-X(1)^2+2X(1)X(6)-X(2)^2+2X(2)X(5)-X(3)^2+2X(3)X(4)-X(4)^2-X(5)^2-X(6)^2} \\
\lambda_2,\lambda_3 &= X(0)+X(7) \\
&\quad - i\sqrt{-X(1)^2+2X(1)X(6)-X(2)^2+2X(2)X(5)-X(3)^2+2X(3)X(4)-X(4)^2-X(5)^2-X(6)^2} \\
\lambda_4,\lambda_5 &= X(0)-X(7) \\
&\quad + i\sqrt{-X(1)^2-2X(1)X(6)-X(2)^2-2X(2)X(5)-X(3)^2-2X(3)X(4)-X(4)^2-X(5)^2-X(6)^2} \\
\lambda_6,\lambda_7 &= X(0)-X(7) \\
&\quad - i\sqrt{-X(1)^2-2X(1)X(6)-X(2)^2-2X(2)X(5)-X(3)^2-2X(3)X(4)-X(4)^2-X(5)^2-X(6)^2}
\end{aligned}
\qquad (32)
$$

and their magnitude is:

$$
\begin{aligned}
|\lambda_0|,|\lambda_1|,|\lambda_2|,|\lambda_3| &= \sqrt{\mathbf{X}_r^T\mathbf{X}_r + \mathbf{X}_r^T\mathbf{J}\mathbf{X}_r} = \sqrt{\sum_{i=0}^{7} X^2(i) - \sum_{i=1}^{6} X(i)X(7-i) + 2X(0)X(7)} = \sqrt{a+b} \\
|\lambda_4|,|\lambda_5|,|\lambda_6|,|\lambda_7| &= \sqrt{\mathbf{X}_r^T\mathbf{X}_r - \mathbf{X}_r^T\mathbf{J}\mathbf{X}_r} = \sqrt{\sum_{i=0}^{7} X^2(i) + \sum_{i=1}^{6} X(i)X(7-i) - 2X(0)X(7)} = \sqrt{a-b}
\end{aligned}
\qquad (33).
$$

The magnitude of the eigenvalues in (33) is therefore the same as the two seminorms defined in (22) and it is non-negative with $\sum_{i=0}^{7} X^2(i) \geq \left| -\sum_{i=1}^{6} X(i)X(7-i) + 2X(0)X(7) \right|$. The cases with zero eigenvalues (and singular $\mathbf{M}_\mathbf{X}$ or $\mathbf{P}_\mathbf{X}$) only occur if one of the two seminorms in (22) is zero, and the conditions required for this to happen have been described in a previous section. Therefore, since all octonion-like numbers with non-zero values for the two seminorms in (22) have an inverse, the octonion-like algebra is a seminormed division algebra, i.e., a division by any number is possible as long as the two seminorms are non-zero.

The following property also applies for the inverse of octonion-like numbers:

$$(\mathbf{XY})^{-1} = \mathbf{Y}^{-1}\mathbf{X}^{-1} \qquad (34),$$

which is easily verified from:

$$\mathbf{XY}(\mathbf{Y}^{-1}\mathbf{X}^{-1})\mathbf{XYY}^{-1}\mathbf{X}^{-1} = 1 \qquad (35)$$

and therefore $(\mathbf{XY})^{-1} = \mathbf{Y}^{-1}\mathbf{X}^{-1}$.

## V. NORMALIZATION OF OCTONION-LIKE NUMBERS AND INVOLUTIONS

We are interested to find an octonion-like number $\mathbf{X}_d$ such that its inverse can be used for a commutative multiplicative normalization:

$$(\mathbf{XX}_d^{-1})(\mathbf{XX}_d^{-1})^* = (\mathbf{X}_d^{-1}\mathbf{X})(\mathbf{X}_d^{-1}\mathbf{X})^* = \mathbf{X}_n\mathbf{X}_n^* = 1 \qquad (36).$$

With the normalised number $\mathbf{X}_n$, the inverse number becomes the conjugate:

$$\mathbf{X}_{n,i} = \mathbf{X}_n^{-1} = \mathbf{X}_n^* \qquad (37).$$

and $\left\|\mathbf{X}_n\right\| = \left\|\mathbf{X}_n^*\right\| = 1$ for the seminorms in (22). Using the previously introduced inverse number, i.e., $\mathbf{X}_{d,r}^{-1} = \mathbf{X}_{i,r}^{-1} = \mathbf{M}_{\mathbf{X}}^{-1}\begin{bmatrix}1\ 0\ 0\ 0\ 0\ 0\ 0\ 0\end{bmatrix}^T$ or $\mathbf{X}_{d,r}^{-1} = \mathbf{X}_{i,r}^{-1} = \mathbf{P}_{\mathbf{X}}^{-1}\begin{bmatrix}1\ 0\ 0\ 0\ 0\ 0\ 0\ 0\end{bmatrix}^T$, would be an obvious solution for this problem, but this would lead to the same normalized value $\mathbf{X}_n = 1\, e_0 = 1$ for any number $\mathbf{X}$ and would be of limited use. Another solution for $\mathbf{X}_d^{-1}$ is described next, which will lead to more interesting involution results.

First we note from the definition of $\mathbf{Z} = \mathbf{XY}$ in (9) that if either $\mathbf{X}$ or $\mathbf{Y}$ has non-zero components only for the real (scalar, non-imaginary) dimensions $e_0, e_7$, then $\mathbf{Z} = \mathbf{XY}$ becomes commutative. We also note that for an octonion-like number $\mathbf{X}$ with non-zero components only for dimensions $e_0, e_7$, $\mathbf{M}_{\mathbf{X}}$ in (7) becomes the sum of a diagonal and an anti-diagonal matrix, therefore $\mathbf{M}_{\mathbf{X}}^{-1}$ also becomes the sum of a diagonal and an anti-diagonal matrix, and consequently the inverse number computed with $\mathbf{X}_{i,r}^{-1} = \mathbf{M}_{\mathbf{X}}^{-1}\begin{bmatrix}1\ 0\ 0\ 0\ 0\ 0\ 0\ 0\end{bmatrix}^T$ is non-zero only for components $e_0, e_7$. Therefore, a simple way to ensure that a product with $\mathbf{X}_d^{-1}$ is commutative is to constrain $\mathbf{X}_d$ to be non-zero only for the $e_0, e_7$ components, i.e., $\mathbf{X}_d = (Ae_0 + Be_7)$. This leads to:

$$(\mathbf{XX}_d^{-1})(\mathbf{XX}_d^{-1})^* = \mathbf{XX}_d^{-1}(\mathbf{X}_d^{-1})^*\mathbf{X}^* = \mathbf{XX}_d^{-1}\mathbf{X}_d^{-1}\mathbf{X}^* = \mathbf{XX}^*\mathbf{X}_d^{-1}\mathbf{X}_d^{-1} = 1 \qquad (38)$$

$$\mathbf{X}_d\mathbf{X}_d = \mathbf{X}\mathbf{X}^* = \sum_{i=0}^{7} X^2(i)\ e_0 + \left( -2\sum_{i=1}^{3} X(i)X(7-i) + 2X(0)X(7) \right) e_7 \tag{39}$$

$$(Ae_0 + Be_7)(Ae_0 + Be_7) = \sum_{i=0}^{7} X^2(i)\ e_0 + \left( -2\sum_{i=1}^{3} X(i)X(7-i) + 2X(0)X(7) \right) e_7 \tag{40}$$

$$A^2 + B^2 = \sum_{i=0}^{7} X^2(i) \tag{41}$$

$$2AB = -2\sum_{i=1}^{3} X(i)X(7-i) + 2X(0)X(7) \qquad (42).$$

Isolating $A$ from (42) and inserting in (41):

$$4B^4 - 4B^2\sum_{i=0}^{7} X^2(i) + \left( -2\sum_{i=1}^{3} X(i)X(7-i) + 2X(0)X(7) \right)^2 = 0 \tag{43}$$

$$4B'^2 - 4B'\sum_{i=0}^{7} X^2(i) + \left( -2\sum_{i=1}^{3} X(i)X(7-i) + 2X(0)X(7) \right)^2 = 0 \tag{44}$$

with $B = \pm\sqrt{B'}$ (both roots are possible and lead to a correct $\mathbf{X}_d$ solution, i.e., $\mathbf{X}_d$ is not unique). Solving for $B'$:

$$B' = \frac{4\sum_{i=0}^{7} X^2(i) + \sqrt{16\left( \sum_{i=0}^{7} X^2(i) \right)^2 - 16\left( -2\sum_{i=1}^{3} X(i)X(7-i) + 2X(0)X(7) \right)^2}}{8} \tag{45}$$

where we use the positive root to guarantee that $B'$ is always positive (the negative root only guarantees non-negativity). Solving for $B', B, A$, we find a commutative number $\mathbf{X}_d = Ae_0 + Be_7$ for which $\mathbf{X}_n\mathbf{X}_n^* = (\mathbf{X}\mathbf{X}_d^{-1})(\mathbf{X}\mathbf{X}_d^{-1})^* = (\mathbf{X}_d^{-1}\mathbf{X})(\mathbf{X}_d^{-1}\mathbf{X})^* = 1$.

For octonion-like numbers $\hat{\mathbf{X}}$ whose non-zero components are only imaginary components $e_1, e_2, e_3, e_4, e_5, e_6$ ("pure" numbers, i.e., real components are zero: $X(0) = 0, X(7) = 0$ ), since $\hat{\mathbf{X}}^* = -\hat{\mathbf{X}}$ the normalization becomes:

$$\left(\hat{\mathbf{X}}_d^{-1}\hat{\mathbf{X}}\right)\left(\hat{\mathbf{X}}_d^{-1}\hat{\mathbf{X}}\right)^* = \hat{\mathbf{X}}_n\hat{\mathbf{X}}_n^* = \left(\hat{\mathbf{X}}_d^{-1}\hat{\mathbf{X}}\right)\left(\hat{\mathbf{X}}^*\hat{\mathbf{X}}_d^{-1}\right) = -\left(\hat{\mathbf{X}}_d^{-1}\hat{\mathbf{X}}\right)\left(\hat{\mathbf{X}}_d^{-1}\hat{\mathbf{X}}\right) = -\hat{\mathbf{X}}_n\hat{\mathbf{X}}_n = 1\ e_0 = 1 \tag{46}$$

$$\hat{\mathbf{X}}_n\hat{\mathbf{X}}_n = -1\ e_0 = -1 \qquad (47).$$

This result may be useful for the development of polar-like representations. Like $\hat{\mathbf{X}}$, $\hat{\mathbf{X}}_n$ also has only imaginary components $e_1, e_2, e_3, e_4, e_5, e_6$ with non-zero values, as the product of a commutative octonion-like number $\hat{\mathbf{X}}_d^{-1}$ (having non-zero coefficients only in dimensions $e_0, e_7$) with a number $\hat{\mathbf{X}}$ (having non-zero coefficients only in dimensions $e_1, e_2, e_3, e_4, e_5, e_6$) also leads to a number $\hat{\mathbf{X}}_n$ with non-zero coefficients only in dimensions $e_1, e_2, e_3, e_4, e_5, e_6$. This can be verified from the definition of $\mathbf{Z} = \mathbf{XY}$ in (9).

We define next the following operation with the same form as rotations in quaternion and octonion algebras:

$$\mathbf{X}^{\mathbf{U}} = \mathbf{U}\mathbf{X}\mathbf{U}^{-1} \tag{48}$$

If $\hat{\mathbf{U}}_n$ is an arbitrary imaginary and "unit" octonion-like number, i.e., it has been normalized such that $\hat{\mathbf{U}}_n\hat{\mathbf{U}}_n = -1$, (48) becomes:

$$\mathbf{X}^{\hat{\mathbf{U}}_n} = \hat{\mathbf{U}}_n\mathbf{X}\hat{\mathbf{U}}_n^{-1} = \hat{\mathbf{U}}_n\mathbf{X}\hat{\mathbf{U}}_n^{*} = -\hat{\mathbf{U}}_n\mathbf{X}\hat{\mathbf{U}}_n \tag{49}$$

which is an involution:

$$-\hat{\mathbf{U}}_n(-\hat{\mathbf{U}}_n\mathbf{X}\hat{\mathbf{U}}_n)\hat{\mathbf{U}}_n = \hat{\mathbf{U}}_n\hat{\mathbf{U}}_n\mathbf{X}\hat{\mathbf{U}}_n\hat{\mathbf{U}}_n = -1 \times \mathbf{X} \times -1 = \mathbf{X} \tag{50}$$

The unit $e_1, e_2, e_3, e_4, e_5, e_6$ are specific cases of such imaginary and "unit" numbers $\hat{\mathbf{U}}_n$. Although not imaginary, the unit $e_0, e_7$ numbers also lead to involutions for $\mathbf{X}^{\mathbf{U}_n} = -\mathbf{U}_n\mathbf{X}\mathbf{U}_n$ or $\mathbf{X}^{\mathbf{U}_n} = \mathbf{U}_n\mathbf{X}\mathbf{U}_n^{-1} = \mathbf{U}_n\mathbf{X}\mathbf{U}_n^{*}$. But in general, normalized non-imaginary numbers $\mathbf{U}_n$ aren't involutions under $\mathbf{X}^{\mathbf{U}_n} = \mathbf{U}_n\mathbf{X}\mathbf{U}_n^{-1} = \mathbf{U}_n\mathbf{X}\mathbf{U}_n^{*}$.

For the definition in (48) we can also derive the following results and properties:

$$\mathbf{X}^{\mathbf{U}} = \mathbf{X}^{\mathbf{U}/a} \quad a \text{ real} \tag{51}$$

$$\mathbf{X}^{\mathbf{U}} = \mathbf{X}^{\mathbf{U}\mathbf{U}_d^{-1}} = \mathbf{X}^{\mathbf{U}_d^{-1}\mathbf{U}} = \mathbf{X}^{\mathbf{U}_n} \tag{52}$$

$$(\mathbf{XY})^{\mathbf{U}} = \mathbf{X}^{\mathbf{U}}\mathbf{Y}^{\mathbf{U}} \tag{53}$$

$$\mathbf{XY} = \mathbf{Y}^{\mathbf{X}}\mathbf{X} = \mathbf{Y}\mathbf{X}^{\mathbf{Y}^*} \tag{54}$$

$$\mathbf{X}^{\mathbf{YZ}} = (\mathbf{X}^{\mathbf{Z}})^{\mathbf{Y}} \tag{55}$$

$$\mathbf{X}^{\mathbf{U}^*} = (\mathbf{X}^*)^{\mathbf{U}} = (\mathbf{X}^{\mathbf{U}})^* \tag{56}$$

For the specific case of the unit octonion-like numbers $e_0(=1), e_1, e_2, e_3, e_4, e_5, e_6, e_7$, we also have:

$$\mathbf{X}^{e_i} = \mathbf{X}^{e_{7-i}} \quad 0 \le i \le 3 \tag{57}$$

$$\begin{aligned}&\mathbf{X}^{e_0} + \mathbf{X}^{e_1} + \mathbf{X}^{e_2} + \mathbf{X}^{e_3} + \mathbf{X}^{e_4} + \mathbf{X}^{e_5} + \mathbf{X}^{e_6} + \mathbf{X}^{e_7} \\ &\quad = 2\mathbf{X}^{e_0} + \mathbf{2X}^{e_1} + 2\mathbf{X}^{e_2} + 2\mathbf{X}^{e_3} = 8X(0)e_0 + 8X(7)e_7\end{aligned} \tag{58}$$

$$(-\mathbf{X}^{e_0} + \mathbf{X}^{e_1} + \mathbf{X}^{e_2} + \mathbf{X}^{e_3} + \mathbf{X}^{e_4} + \mathbf{X}^{e_5} + \mathbf{X}^{e_6} - \mathbf{X}^{e_7})^* = (-2\mathbf{X}^{e_0} + 2\mathbf{X}^{e_1} + 2\mathbf{X}^{e_2} + 2\mathbf{X}^{e_3})^* = 4\mathbf{X} \quad (59).$$

Finally, using (48) with non-zero $\mathbf{U}$ generates an alternative basis:

$$\begin{aligned}\mathbf{X}^{\mathbf{U}} &= \mathbf{U}\mathbf{X}\mathbf{U}^{-1} = \mathbf{U}(x_0 + x_1e_1 + x_2e_2 + x_3e_3 + x_4e_4 + x_5e_5 + x_6e_6 + x_7e_7)\mathbf{U}^{-1} \\ &= x_0 + x_1\mathbf{U}e_1\mathbf{U}^{-1} + x_2\mathbf{U}e_2\mathbf{U}^{-1} + x_3\mathbf{U}e_3\mathbf{U}^{-1} + x_4\mathbf{U}e_4\mathbf{U}^{-1} + x_5\mathbf{U}e_5\mathbf{U}^{-1} + x_6\mathbf{U}e_6\mathbf{U}^{-1} + x_7\mathbf{U}e_7\mathbf{U}^{-1} \\ &= x_0 + x_1e_1^{\mathbf{U}} + x_2e_2^{\mathbf{U}} + x_3e_3^{\mathbf{U}} + x_4e_4^{\mathbf{U}} + x_5e_5^{\mathbf{U}} + x_6e_6^{\mathbf{U}} + x_7e_7^{\mathbf{U}}\end{aligned} \tag{60}$$

where the resulting $e_i^{\mathbf{U}}$ terms still obey the rules of the octonion-like algebra:

$$e_0^{\mathbf{U}} = e_0 = 1 \qquad e_7^{\mathbf{U}} = e_7 \qquad e_i^{\mathbf{U}}e_i^{\mathbf{U}} = -1 \;\; 1 \le i \le 6 \qquad e_1^{\mathbf{U}}e_2^{\mathbf{U}}e_3^{\mathbf{U}}e_4^{\mathbf{U}}e_5^{\mathbf{U}}e_6^{\mathbf{U}} = 1 \qquad \left\| e_i^{\mathbf{U}} \right\| = 1 \;\; 0 \le i \le 7 \quad (61).$$

## VI. CONCLUSION

Using elementary linear algebra, the associative octonion-like algebra (split-biquaternion algebra, even subalgebra of Clifford algebra $Cl_{4,0}(\mathbb{R})$) was shown to be a seminormed composition algebra over $\mathbb{R}$ and a seminormed division algebra over $\mathbb{R}$ (i.e., division by any number is possible as long as two seminorms are non-zero). This does not contradict Hurwitz's theorem. It is in general a non-commutative algebra (just like quaternions and octonions), although some operations such as number inverse, normalization or multiplication with numbers having no imaginary components are commutative. This paper also provided results for the computation and existence of inverse numbers in the octonion-like algebra, as well as for number normalization and some involutions. The results were developed using elementary linear algebra, which can facilitate software implementations.

## VII. ACKNOWLEDGEMENTS

The authors would like to thank Dr. Joy Christian and Dr. Richard D. Gill for helpful comments and clarifications on different versions of this paper.

## VIII. REFERENCES

[1] J. Christian, "Eight-dimensional Octonion-like but Associative Normed Division Algebra," *arXiv:1908.06172v8 [math.GM]* , November 2020.

[2] J. Christian, "Quantum Correlations Are Weaved by the Spinors of the Euclidean Primitives," *R. Soc. Open Sci., ,* Vols. 5, 180526 , 2018.

[3] "Hurwitz's theorem (composition algebras) (Wikipedia)," [Online]. Available: https://en.wikipedia.org/wiki/Hurwitz%27s_theorem_(composition_algebras).

[4] R. D. Gill, "Does Geometric Algebra Provide a Loophole to Bell's Theorem?," *Entropy,* Vols. 22, 61; doi:10.3390/e22010061, 2020.

[5] "Split-biquaternion (Wikipedia)," [Online]. Available: https://en.wikipedia.org/wiki/Split-biquaternion.

[6] A. N. Lasenby, "A 1d Up Approach to Conformal Geometric Algebra: Applications in Line Fitting and Quantum Mechanics," *Adv. Appl. Clifford Algebras* , vol. 30:22, 2020.

## IX. APPENDIX: PROOF FOR PRESERVATION OF SEMINORMS PRODUCT UNDER OCTONION-LIKE NUMBER MULTIPLICATION

We need to show that:

$$\left\|\mathbf{Z}\right\|^2 = \left\|\mathbf{XY}\right\|^2 = \left\|\mathbf{X}\right\|^2 \left\|\mathbf{Y}\right\|^2 \qquad (62),$$

where $\mathbf{X},\mathbf{Y},\mathbf{Z}$ are octonion-like numbers as before. With the two seminorms for octonion-like numbers in (22), this is equivalent to:

$$\begin{aligned}&\sum_{i=0}^{7} Z^2(i) - 2\lambda\sum_{i=1}^{3} Z(i)Z(7-i) + 2\lambda Z(0)Z(7)\\&=\left(\sum_{i=0}^{7} X^2(i) - 2\lambda\sum_{i=1}^{3} X(i)X(7-i) + 2\lambda X(0)X(7)\right)\left(\sum_{i=0}^{7} Y^2(i) - 2\lambda\sum_{i=1}^{3} Y(i)Y(7-i) + 2\lambda Y(0)Y(7)\right)\end{aligned} \quad (63),$$

with $\lambda = \pm 1$. In matrix form, since $\left\|\mathbf{X}\right\|^2 = \mathbf{X}_r^T\mathbf{X}_r + \lambda\mathbf{X}_r^T(\mathbf{J}\mathbf{X}_r)$ we have:

$$\begin{aligned}\mathbf{Z}_r^T\mathbf{Z}_r + \lambda\mathbf{Z}_r^T(\mathbf{J}\mathbf{Z}_r) &= \left(\mathbf{X}_r^T\mathbf{X}_r + \lambda\mathbf{X}_r^T(\mathbf{J}\mathbf{X}_r)\right)\left(\mathbf{Y}_r^T\mathbf{Y}_r + \lambda\mathbf{Y}_r^T(\mathbf{J}\mathbf{Y}_r)\right)\\&= \left(\mathbf{Y}_r^T\mathbf{Y}_r + \lambda\mathbf{Y}_r^T(\mathbf{J}\mathbf{Y}_r)\right)\left(\mathbf{X}_r^T\mathbf{X}_r + \lambda\mathbf{X}_r^T(\mathbf{J}\mathbf{X}_r)\right)\end{aligned} \quad (64).$$

Since $\mathbf{Z}_r = \mathbf{M}_{\mathbf{X}}\mathbf{Y}_r$, we can write:

$$\begin{aligned}&\mathbf{Y}_r^T\mathbf{M}_{\mathbf{X}}^T\mathbf{M}_{\mathbf{X}}\mathbf{Y}_r + \lambda\mathbf{Y}_r^T\mathbf{M}_{\mathbf{X}}^T(\mathbf{J}\mathbf{M}_{\mathbf{X}}\mathbf{Y}_r)\\&= \mathbf{Y}_r^T\mathbf{Y}_r\mathbf{X}_r^T\mathbf{X}_r + \lambda\mathbf{Y}_r^T\mathbf{Y}_r\mathbf{X}_r^T(\mathbf{J}\mathbf{X}_r) + \lambda\mathbf{Y}_r^T(\mathbf{J}\mathbf{Y}_r)\mathbf{X}_r^T\mathbf{X}_r + \mathbf{Y}_r^T(\mathbf{J}\mathbf{Y}_r)\mathbf{X}_r^T(\mathbf{J}\mathbf{X}_r)\end{aligned} \quad (65)$$

and

$$\begin{aligned}&\mathbf{M}_{\mathbf{X}}^T\mathbf{M}_{\mathbf{X}}\mathbf{Y}_r + \lambda\mathbf{M}_{\mathbf{X}}^T\mathbf{J}\mathbf{M}_{\mathbf{X}}\mathbf{Y}_r\\&= \mathbf{Y}_r(\mathbf{X}_r^T\mathbf{X}_r) + \lambda\mathbf{Y}_r(\mathbf{X}_r^T\mathbf{J}\mathbf{X}_r) + \lambda\mathbf{J}\mathbf{Y}_r(\mathbf{X}_r^T\mathbf{X}_r) + \mathbf{J}\mathbf{Y}_r(\mathbf{X}_r^T\mathbf{J}\mathbf{X}_r)\end{aligned} \quad (66).$$

This condition holds true if the two following conditions are true:

$$\mathbf{M}_{\mathbf{X}}^T\mathbf{M}_{\mathbf{X}}\mathbf{Y}_r = \mathbf{Y}_r(\mathbf{X}_r^T\mathbf{X}_r) + \mathbf{J}\mathbf{Y}_r(\mathbf{X}_r^T\mathbf{J}\mathbf{X}_r) \quad (67)$$

$$\lambda\mathbf{M}_{\mathbf{X}}^T\mathbf{J}\mathbf{M}_{\mathbf{X}}\mathbf{Y}_r = \lambda\mathbf{Y}_r(\mathbf{X}_r^T\mathbf{J}\mathbf{X}_r) + \lambda\mathbf{J}\mathbf{Y}_r(\mathbf{X}_r^T\mathbf{X}_r) \quad (68).$$

The 1st condition is easily shown using the previous result $\mathbf{M}_{\mathbf{X}}^T\mathbf{M}_{\mathbf{X}} = (\mathbf{X}_r^T\mathbf{X}_r)\mathbf{I} + (\mathbf{X}_r^T(\mathbf{J}\mathbf{X}_r))\mathbf{J}$:

$$\begin{aligned}\mathbf{M}_{\mathbf{X}}^T\mathbf{M}_{\mathbf{X}}\mathbf{Y}_r &= ((\mathbf{X}_r^T\mathbf{X}_r)\mathbf{I} + (\mathbf{X}_r^T(\mathbf{J}\mathbf{X}_r))\mathbf{J})\mathbf{Y}_r\\&= (\mathbf{X}_r^T\mathbf{X}_r)\mathbf{Y}_r + (\mathbf{X}_r^T\mathbf{J}\mathbf{X}_r)\mathbf{J}\mathbf{Y}_r\\&= \mathbf{Y}_r(\mathbf{X}_r^T\mathbf{X}_r) + \mathbf{J}\mathbf{Y}_r(\mathbf{X}_r^T\mathbf{J}\mathbf{X}_r)\end{aligned} \quad (69).$$

For the 2nd condition, first we show that $\mathbf{M}_{\mathbf{X}}^T\mathbf{J}\mathbf{M}_{\mathbf{X}} = \mathbf{J}\mathbf{M}_{\mathbf{X}}^T\mathbf{M}_{\mathbf{X}}$, making use of $\mathbf{J}\mathbf{J} = \mathbf{I}$:

$$\mathbf{M}_{\mathbf{X}}^T\mathbf{J}\mathbf{M}_{\mathbf{X}}\mathbf{Y}_r = \mathbf{Y}_r(\mathbf{X}_r^T\mathbf{J}\mathbf{X}_r) + \mathbf{J}\mathbf{Y}_r(\mathbf{X}_r^T\mathbf{X}_r) \quad (70)$$

$$\begin{aligned}\mathbf{J}\mathbf{M}_{\mathbf{X}}^{T}\mathbf{J}\mathbf{M}_{\mathbf{X}}\mathbf{Y}_{r} &= \mathbf{J}\mathbf{Y}_{r}\left(\mathbf{X}_{r}^{T}\mathbf{J}\mathbf{X}_{r}\right)+\mathbf{J}\mathbf{J}\mathbf{Y}_{r}\left(\mathbf{X}_{r}^{T}\mathbf{X}_{r}\right)\\ &= \mathbf{J}\mathbf{Y}_{r}\left(\mathbf{X}_{r}^{T}\mathbf{J}\mathbf{X}_{r}\right)+\mathbf{Y}_{r}\left(\mathbf{X}_{r}^{T}\mathbf{X}_{r}\right)\\ &= \mathbf{M}_{\mathbf{X}}^{T}\mathbf{M}_{\mathbf{X}}\mathbf{Y}_{r}\end{aligned} \tag{71}$$

$$\mathbf{J}\mathbf{M}_{\mathbf{X}}^{T}\mathbf{J}\mathbf{M}_{\mathbf{X}} = \mathbf{M}_{\mathbf{X}}^{T}\mathbf{M}_{\mathbf{X}} \tag{72}$$

$$\mathbf{J}\mathbf{J}\mathbf{M}_{\mathbf{X}}^{T}\mathbf{J}\mathbf{M}_{\mathbf{X}} = \mathbf{J}\mathbf{M}_{\mathbf{X}}^{T}\mathbf{M}_{\mathbf{X}} \tag{73}$$

$$\mathbf{M}_{\mathbf{X}}^{T}\mathbf{J}\mathbf{M}_{\mathbf{X}} = \mathbf{J}\mathbf{M}_{\mathbf{X}}^{T}\mathbf{M}_{\mathbf{X}} \tag{74.}$$

Then the 2nd condition is easily verified, again using $\mathbf{M}_{\mathbf{X}}^{T}\mathbf{M}_{\mathbf{X}} = (\mathbf{X}_{r}^{T}\mathbf{X}_{r})\mathbf{I}+(\mathbf{X}_{r}^{T}(\mathbf{J}\mathbf{X}_{r}))\mathbf{J}$:

$$\begin{aligned}\mathbf{M}_{\mathbf{X}}^{T}\mathbf{J}\mathbf{M}_{\mathbf{X}}\mathbf{Y}_{r} &= \mathbf{J}\mathbf{M}_{\mathbf{X}}^{T}\mathbf{M}_{\mathbf{X}}\mathbf{Y}_{r} = \mathbf{J}\Big((\mathbf{X}_{r}^{T}\mathbf{X}_{r})\mathbf{I}+(\mathbf{X}_{r}^{T}(\mathbf{J}\mathbf{X}_{r}))\mathbf{J}\Big)\mathbf{Y}_{r}\\ &= \mathbf{J}\mathbf{Y}_{r}(\mathbf{X}_{r}^{T}\mathbf{X}_{r})+\mathbf{Y}_{r}(\mathbf{X}_{r}^{T}\mathbf{J}\mathbf{X}_{r})\end{aligned} \tag{75.}$$